\documentstyle[12pt]{article}

\input amssym.def
\input amssym.tex

\newcommand{\nc}{\newcommand}

\newcommand{\Tr}{\,{\rm Tr}\,}

\newcommand{\Hom}{\,{\rm Hom}\,}

\newcommand{\Ker}{ \,{\rm Ker} \,}

\newcommand{\bla}{\phantom{bbbbb}}
\newcommand{\onebl}{\phantom{a} }

\newcommand{\half}{ {\frac{1}{2} } }
\newcommand{\vol}{ \,{\rm vol}\, }

\newcommand{\beq}{\begin{equation}}
\newcommand{\eeq}{\end{equation}}
\newcommand{\beqst}{\begin{equation*}}
\newcommand{\eeqst}{\end{equation*}}
\newcommand{\barr}{\begin{array}}
\newcommand{\earr}{\end{array}}
\newcommand{\beqar}{\begin{eqnarray}}
\newcommand{\eeqar}{\end{eqnarray}}
\newtheorem{theorem}{Theorem}[section]

\newtheorem{corollary}[theorem]{Corollary}

\newtheorem{lemma}[theorem]{Lemma}
\newtheorem{prop}[theorem]{Proposition}

\newtheorem{remit}[theorem]{Remark}

\newcommand{\RR}{{\Bbb R }}
\newcommand{\CC}{{\Bbb C }}
\newcommand{\ZZ}{{\Bbb Z }}

\newcommand{\cala}{{\mbox{$\cal A$}}}

\newcommand{\calf}{{\mbox{$\cal F$}}}

\newcommand{\calh}{{\mbox{$\cal H$}}}
\newcommand{\cali}{{\mbox{$\cal I$}}}

\newcommand{\call}{{\mbox{$\cal L$}}}

\newcommand{\cals}{{\mbox{$\cal S$}}}

\def\b{\beta}
\def\g{\gamma}

\def\e{\epsilon}

\def\k{\kappa}
\def\l{\lambda}

\def\G{\Gamma}

\def\L{\Lambda}

\nc{\Proof}{\noindent{\em Proof:} }
\nc{\Lie}{ {\rm Lie} }

\nc{\liek}{{ \bf k} }
\nc{\lieks}{\liek^*}
\nc{\symm}{M}
\nc{\hk}{H^*_K}

\nc{\quant}{\calh}
\nc{\quantr}{\calh_{\rm red} }
\nc{\lred}{\call_{\rm red} }
\nc{\mred}{\symm_{\rm red} }
\nc{\td}{{\rm Td} }
\nc{\tdn}{{\rm Td}^{<n> } }
\nc{\tdnk}{{\rm Td}_K^{<n> } }
\nc{\tdnt}{{\rm Td}_T^{<n> } }
\nc{\tdnki}{ ({\rm Td}_K^{-1 })^{<n> }  }
\nc{\tdnti}{ ({\rm Td}_T^{-1 })^{<n> }  }
\nc{\tdk}{{\rm Td}_K }
\nc{\tdt}{{\rm Td}_T }

\nc{\tdki}{({\rm Td}_K)^{-1}}
\nc{\tdti}{({\rm Td}_T)^{-1}}

\nc{\lieksa}{\lieks_{{\rm ad}} }
\nc{\lieka}{\liek_{{\rm ad}} }
\nc{\isq}{ { i } }
\nc{\isqrm}{~}
\nc{\rrk}{RR^K}
\nc{\rrt}{RR^T}
\nc{\rr}{RR}
\nc{\chk}{{\rm ch}_K}
\nc{\cht}{{\rm ch}_T}
\nc{\ch}{ {\rm ch} }

\nc{\proj}{\pi}
\nc{\zloc}{\mu^{-1}(0)}
\nc{\nf}{\nu_F}
\nc{\lx}{\l(X) }
\nc{\rlx}{\frac{d \l(X)}{\l(X) }  }
\nc{\srlx}{ {d \l(X)} / {\l(X) }  }
\nc{\intlat} {\L^I}

\nc{\hg}{ {\hat{\g} } }
\nc{\chtd}{ { \int_{\mred} \ch(\lred) \td(\mred) } }
\nc{\sumpl}{ { \sum_{F \in \calf_+} } }

\nc{\nplf}{{n_{F,+} } }
\nc{\nminf}{{n_{F,-} } }
\nc{\liet}{{\bf t}}
\nc{\nusym}{\varpi}
\nc{\rfe}{r_F^{\eta} }
\nc{\liets}{\liet^*}
\nc{\expmf}{e^{\isq \mu_T(F)(X) } }
\nc{\expbfj}{e^{\isq \bfj (X) } }
\nc{\expmfs}{e^{\isq \mu_T(F)X } }
\nc{\sumf}{\sum_{F \in \calf_+} }
\nc{\mf}{\mu_T(F)}
\nc{\tck}{{C^K_1}}
\nc{\cko}{C^K_2}
\nc{\bom}{{\omega_\liek}}
\nc{\fk}{F_K}
\nc{\ft}{F_T}

\nc{\sss}{s}
\nc{\wf} { { { \cal W}_F } }

\nc{\charr}{\chi}
\nc{\bfj}{ {\beta_{F,j} } }
\nc{\nfj}{\nu_{F,j} }
\nc{\hatc}{ \hat{\CC} }
\nc{\res}{{\rm res} }

\nc{\ms}{M_s}
\nc{\msoa}{{ \ms^{0,a} } }
\nc{\msa}{{ \ms^{a} } }
\nc{\msoab}{{    M^a_{s,{\rm red}}     } }
\nc{\isa}{\cali^{s,a} }
\nc{\kappind}{\kappa_a}
\nc{\mua}{\mu_a}
\nc{\betak}{ {\b_{s,a,k} } }
\nc{\cnus}{ { c_1(\nu_{s,a,k} ) } }
\nc{\cnust}{ {  \cnus_T } }
\nc{\nusak}{{\nu_{s,a,k} } }
\nc{\sumfms}{ { \sum_{F \in \calf_+ : F\subset  \msa} } }
\nc{\isat}{\isa_T}
\nc{\stabo}{n_0}

\newcommand{\renorm}{{ \setcounter{equation}{0} }}

\begin{document}

\title{On localization
and Riemann-Roch numbers for symplectic quotients}

\author{Lisa C. Jeffrey \\
Mathematics Department \\ Princeton University
\\ Princeton NJ 08544, USA \thanks{This material is based on work
supported by the National Science Foundation under Grant No.
DMS-9306029.}
\\and \\
Frances C. Kirwan \\ Balliol College, Oxford OX1 3BJ, UK}

\date{September 1994; revised March 1995}

\maketitle
\begin{abstract}
Suppose $(M,\omega)$ is a compact symplectic manifold acted on by a
compact Lie group $K$ in a Hamiltonian fashion, with
moment map $\mu: M \to \Lie(K)^*$ and Marsden-Weinstein
reduction $\mred = \mu^{-1}(0)/K$. In our earlier
paper \cite{JK1}, under the assumption that $0$ is  a
regular value of $\mu$, we proved a formula (the residue
formula) for $\eta_0 e^{\omega_0}[\mred]$ for
any $\eta_0 \in H^*(\mred)$, where $\omega_0$ is
the induced symplectic form on $\mred$. This formula is given
in terms of the restrictions of classes in the equivariant
cohomology $H^*_T(\symm) $ of $\symm$ to the
components of the fixed point set of
a maximal torus $T$ in $\symm$.

In this paper, we assume that
$\symm$ has a $K$-invariant K\"ahler structure.
We  apply the residue formula in the
special case $\eta_0 = \td(\mred)$; when $K$ acts freely on $\mu^{-1}(0)$
this yields a formula for the Riemann-Roch number $\rr (\lred)$
of a holomorphic
line bundle $\lred$ on $\mred$
that descends from a holomorphic
line bundle $\call$ on $\symm$ for which $c_1(\call) = \omega$.
More generally when 0 is a regular value of $\mu$, so that $\mred$
is an orbifold and $\lred$ is an orbifold bundle, Kawasaki's Riemann-Roch
theorem for orbifolds can be applied.
Using the holomorphic Lefschetz formula we similarly
obtain a formula for the  $K$-invariant Riemann-Roch
number  $\rrk(\call) $ of $\call$. In the case when
the maximal torus $T$ of $K$ has dimension one (except
in a few special circumstances), we show
the two formulas are the same. Thus in this special case the
residue formula is equivalent to the result of Guillemin and
Sternberg \cite{gsgq} that $\rr(\lred) = \rrk(\call)$.
\end{abstract}

\renorm
\section{Introduction}

Let $M$ be a compact symplectic manifold of (real) dimension $2m$,
acted on in a Hamiltonian fashion by a compact connected
Lie group $K$
with maximal torus $T$, and let $\liek$ and $\liet$ denote the
Lie algebras of  $K$ and $T$.  Let $\mu: M \to \lieks$
be a moment map for this action. The reduced space $\mred$ is defined
as
$$ \mred = \mu^{-1}(0)/K. $$
We shall  assume throughout
this paper that $0$ is a regular value of $\mu$, so that
$\mred$ is a symplectic orbifold; it
has at worst finite quotient singularities
and the symplectic form $\omega$ on $M$ induces a symplectic
$\omega_0$ on $\mred$.

There is a natural surjective ring homomorphism \cite{Ki1}
$\k_0: \hk(M) \to H^*(\mred)$, where $\hk(M)$ is the
$K$-equivariant cohomology of $M$.
The main result of
\cite{JK1} was the residue formula  (Theorem 8.1),
which  for any $\eta_0 \in H^*(\mred)$
gives a formula for the evaluation of $\eta_0 e^{ \omega_0}$ on the
fundamental class of $\mred$.
This is given in terms of the restrictions $i_F^* \eta $
 to components $F$ of the fixed point
set    of $T$ in $\symm$, for  any class $\eta \in \hk(M)$ which maps to
$\eta_0$ under $\k_0$. The residue formula
is an application of the localization
theorem of Berline and Vergne,  a result on the
 equivariant cohomology
of torus actions \cite{BV1}, for which a
topological proof was later given by
Atiyah and Bott \cite{abmm}. The residue formula is
related to a result of Witten (the nonabelian
localization theorem \cite{tdg}):
like the residue formula, Witten's theorem   expresses $\eta_0
e^{\omega_0}[\mred] $ in terms of
appropriate data on $\symm$.

In this paper
we assume also that
 there exists a line bundle $\call$
on $M$ for which $c_1(\call) = \omega$, with the
action of $K$ on $M$ lifting to an action on the total
space of $\call$ (such a lift exists because the action
of $K$ on $M$ is Hamiltonian; we choose the lift to be compatible with
the chosen moment map).
Under the assumption that $K$ acts freely on
$\mu^{-1} (0)$,
we get a line bundle $\lred$ on $\mred$ whose first Chern
class is $\omega_0$, where $\omega_0$ is the induced symplectic
form on $\mred$.
In the more general case, $\lred$   is only an
orbifold bundle.

Suppose also that there exists a $K$-invariant
K\"ahler structure on $\symm$; more precisely a
complex structure  compatible
with $\omega $ and preserved by
the
action of $K$.
The bundle $\call $ then acquires a holomorphic structure in a
standard manner, and
we  define the {\em quantizations}
$\quant$ and $\quantr$ to be the virtual vector spaces
\beq \label{2.1} \quant = \oplus_{j \ge 0 }
(-1)^j H^j(\symm, \call) \eeq
and
\beq \label{2.1b} \quantr = \oplus_{j \ge 0 }
(-1)^j H^j(\mred, \lred). \eeq
The space $\quant$ is a virtual representation of $K$.
The {\rm Riemann-Roch numbers} $\rrk(\call)$ and
$\rr (\lred) $ are defined by
\beq \label{2.2}
\rrk (\call) = \sum_{j \ge 0 } (-1)^j \dim H^j (\symm, \call)^K \eeq
\beq \label{2.2b}
\rr (\lred) = \sum_{j \ge 0 } (-1)^j \dim H^j (\mred, \lred). \eeq

The main result of this paper is that when the dimension of the maximal
torus $T$ is one, except in a few special circumstances,
a particular case of our residue formula is equivalent to the statement that
these two Riemann-Roch numbers are equal:
\beq \label{2.02} \rr (\lred) = \rrk(\call). \eeq
This statement was proved by Guillemin and Sternberg [14] under some
additional positivity hypotheses on $\call$, and was conjectured by them to
hold more generally. It has been called the quantization conjecture: that
quantization commutes with reduction.

In this paper we show the following:

\noindent{\bf Theorem 6.2:} {\em Suppose
$K$ is a compact connected group of rank one. Let
$K$ act in a Hamiltonian fashion on the symplectic manifold
$\symm$, with a moment map $\mu$ for the action of $K$
such that $0$ is a regular value of $\mu$.
 If $K = SO(3)$, suppose also
 that there exists a component $F$ of the fixed point set
$M^T$ of the maximal torus $T$
such that the constant value
taken by the $T$-moment map $\mu_T$ on $F$satisfies
$|\mf| > 1$, and if  $K=SU(2)$
suppose that there is an $F$ for which $|\mf|>2$ and that
there is no $F$ with $\mf=\pm1=n_{F,\pm}$
where $n_{F,\pm}$ is the sum of the positive (respectively negative)
weights for the action of $T$ on
the normal to $F$ in $M$. Then $\rrk (\call) = \rr (\lred).$  }

 Our original motivation for considering
Riemann-Roch numbers  was to provide a link between the residue
we had defined and more standard definitions
of residues in algebraic geometry (such as the Grothendieck
residue \cite{harts}).
 We had defined the residue as the evaluation at $0$ (suitably interpreted)
of the Fourier transform of a particular function on $\liet$:
 in the case when $T$ has rank one, this may
be identified with the residue at $0$ of a meromorphic
function on $\CC$ whose poles occur only at $0$. Moreover, when $T$
has
rank one, the special case of our residue that arises in
computing $\rr (\lred)$ may be recast as the residue of a meromorphic
1-form on the Riemann sphere $\hatc$ at one of its poles. This same
expression
arises when one computes $\rrk(\call)$ by using the holomorphic
Lefschetz formula to give a formula for the character  $\charr(k)$ of the
action of an element $k$ of $K$ on $\quant$ and then integrating
$\charr(k)$ over the group $K$ to get the dimension of the
$K$-invariant
subspace $\quant^K$.
Since writing this paper we have
found that it is possible to extend its methods to treat
the case when $K$ has higher rank (see \cite{JKQ});
however the arguments become more involved.

Since we first began considering the application of the residue formula
to Riemann-Roch numbers, several papers have appeared which
extend the Guillemin-Sternberg result to a wider class of situations,
and in which the main tool is localization in equivariant cohomology.
There are two approaches, one due to Guillemin \cite{G}, the other
due to Vergne \cite{V}.

Guillemin's proof uses the
residue formula to reduce
the verification of (1.5) to a
combinatorial identity involving counting lattice points in polyhedra.
Guillemin then observes that
this identity is known when $K$ is a torus acting in
a quasi-free manner.
Meinrenken \cite{mein} has subsequently extended this proof to torus actions
which need not be quasi-free.

As has been pointed out by Guillemin (\cite{G}, Section 3), the
application of the residue formula  to yield a formula for
$\rr(\lred)$ requires only that there exist an {\em almost complex}
structure on $\symm$ compatible with
the action of $K$: such a structure enables one to define
a spin-$\CC$  Dirac operator which can be used to define the virtual
vector space $\quant$.  Guillemin and Sternberg's
original proof \cite{gsgq}, on the other hand,
 depends on the existence of a K\"ahler structure on $\symm$ and
on some positivity hypotheses that are not necessary in the approaches
based on equivariant cohomology.
Thus the use of localization in equivariant cohomology
extends  Guillemin and Sternberg's original result
 to a more general situation. The observation that it
suffices to assume the existence of a $K$-invariant compatible
almost complex structure on $\symm$
applies likewise to the proof we shall present.

Vergne [23] has
given a different proof of the Guillemin-Sternberg conjecture when
$K$ is a torus, also using ideas based
on localization in
equivariant cohomology.
Her proof  likewise does not  require positivity
hypotheses or the existence of a
$K$-invariant K\"ahler structure on $\symm$.

In later work, Meinrenken
\cite{mein2} has proved the Guillemin-Sternberg formula for general compact
nonabelian groups $K$; the only hypothesis he imposes is that
the symmetric quadratic form on the tangent space to $\symm$
given by the symplectic
form and the almost complex structure should be a
metric (i.e., that it should be positive definite).

Although many features of the rank one case are quite special, and
although the proofs of Vergne and Guillemin-Meinrenken described
above apply in much greater generality, we felt
nevertheless that it was
instructive to give a written
account of  our approach to this case since it is
simple and self-contained.

The layout of this paper is as follows. In Section 2 we review some basic facts
about equivariant cohomology, and find an equivariant cohomology class $\eta$
on $M$ mapping to the Todd class of $\mred$ under the natural surjection
$\k_0:\hk(M) \to H^*(\mred)$, so that
$$\eta_0 = \td(\mred).$$
By the Riemann-Roch formula we have
\beq \label{3.1} RR(\lred) = \eta_0 e^{\omega_0}[\mred], \eeq
provided $\mred$ is a manifold,
which is true if $K$ acts freely on $\mu^{-1}(0)$;
in the more general
case $RR(\lred)$ is given by Kawasaki's Riemann-Roch theorem for
orbifolds (Theorem 6.1).
In Section 3 we apply the residue formula to give a formula
for the right hand side of (\ref{3.1}
as a sum over the
components of the fixed point
set $M^T$ of $T$.
In Section 4, we apply
the holomorphic Lefschetz
formula to obtain a similar fixed point sum for $RR^K(\lred)$;
finally in Sections 5 and 6 we identify the two expressions.

\noindent{\em Acknowledgement:} We thank Miles Reid, who originally
proposed  that we should investigate the possibility
of trying to reformulate the residue in \cite{JK1}
in more algebro-geometric terms, and suggested
that there should be a relation
between our residue formula and the Riemann-Roch theorem.

\renorm
\section{Preliminaries}

If $M$ is a compact oriented manifold acted on by a compact
connected Lie group $K$, the $K$-equivariant cohomology $\hk(M)$
of $M$ may be identified with the cohomology of the following
chain complex (see Chapter 7 of \cite{BGV}):
\beq \Omega^*_K(M) = \Bigl ( S(\lieks) \otimes \Omega^*(M) \Bigr )^K
\eeq
equipped with the  differential\footnote{This (nonstandard) definition
of the
equivariant cohomology differential is different from that used in
\cite{JK1} but consistent with that used in \cite{tdg}. We have
found it convenient  to introduce this definition to obtain
consistency with  the formulas in Section 4.}

\beq \label{3.001}
D = d - \isq \iota_{X_M}  \eeq
where $X_M$ is the vector field on $M $ generated by the
action of $X \in \liek$. The  natural map
\beq \label{3.002}
\tau_M: \hk(\symm) \to H^*_T (\symm) \eeq
corresponds to the restriction map
$$ \Bigl ( S(\lieks) \otimes \Omega^*(M) \Bigr )^K
\to S(\liets) \otimes \Omega^*(M).  $$

We shall make use of equivariant characteristic classes on
$\symm$: for their definition see Section 7.1 of \cite{BGV}.
\begin{prop} \label{p2.1}
The equivariant Chern character of the line bundle $\call$ is
$$ \chk(\call) (X) =  \Tr (e^{ \omega + \isq (\mu,X) } )  \in
\hat{\Omega}^*_K(\symm). $$
\end{prop}
Here, $X$ is a parameter in
$\liek$, and
$\hat{\Omega}^*_K(\symm) $ is the formal completion of
the space ${\Omega}^*_K(\symm)$  of equivariant differential forms
on $\symm$.

Suppose $F$ is a component of the fixed point set of $T$ in
$\symm$.
We may (formally) decompose the normal
bundle $\nu_F$ to $F$ (using the splitting principle if necessary) as
a sum of line bundles
$\nu_F = \sum_{j = 1}^N \nfj $, in such a way that
$T$ acts on $\nfj$ with weight $\bfj \in \liets$.\footnote{Throughout
this paper we shall use the convention that weights $\bfj \in \liets$
send the integer lattice $\intlat = \Ker(\exp: \liet \to T) $
to $\ZZ$.}
The $T$-equivariant
Euler class $e_F$ of the normal bundle $\nu_F$
is then defined for $X \in \liet$ by
\beq \label{3.021}
 e_F(X)  = \prod_{j  = 1}^N \bigl (c_1 (\nfj) + \isq \bfj(X) \bigr). \eeq

Recall that
the Todd class of a vector bundle $V$  is given in terms of the Chern
roots
$x_l$ by
$$\td  (V) = \prod_l \frac{x_l}{1 - e^{- x_l} }
= \sum_{j \ge 0 } \td_j(V), $$
where $\td_j $ is a  homogeneous polynomial  of degree $j$ in the
$x_l$.
If the    Todd class is given in terms of the Chern
roots by
$$\td  = \tau (x_1, \dots, x_N) $$
then
the $T$-equivariant
  Todd class of the normal bundle  $\nf $ is given for $X \in \liet$
by
\beq \label{4.9} \td_T (\nu_F ) (X) =
\tau \Bigl ( c_1 (\nu_{F,1})  + \isq \beta_{F,1}(X) , \dots,
c_1 (\nu_{F,N})  + \isq \beta_{F,N}(X)  \Bigr )  . \eeq

We may also define the $K$-equivariant Todd class $\tdk(V)$
of any $K$-equivariant vector bundle $V$ on $M$, and in particular
the equivariant Todd class
$\tdk (\symm) $
$= \tdk(T\symm)$  of $M$.  We have
$ \tau_M (\tdk(V) ) = \tdt(V)$  and $\tau_M
(\chk(\call) ) = \ch_T (\call) $, where
$\tau_M$ is the natural  map introduced at (\ref{3.002}).
Moreover one may define the inverse equivariant
Todd class $$\tdki(V) = \sum_{j = 0 }^\infty
(\tdki)_j (V) $$
 as the  equivariant extension of  the class $\td^{-1}(V)$ given in
terms of the Chern roots by
$$\td^{-1}(V) =
\prod_l \frac{1 - e^{- x_l} } {x_l}.
 $$

The surjective ring homomorphism
$\k_0:\hk(\symm) \to H^*(\mred)$ mentioned in the
introduction is the composition of
the restriction map from $\hk(M)$ to $\hk(\mu^{-1}(0))$
and the natural isomorphism
from $\hk(\mu^{-1}(0))$ to
$H^*(\mred)$ which exists since
$K$ acts locally freely on
$\mu^{-1}(0)$ and we are working with cohomology with complex
coefficients. This surjection
is zero on $H^j_K(M)$ for any $j> dim_{\bf R}(\mred)$, and
so it makes sense to apply
$\k_0$ to formal equivariant cohomology classes such as the
equivariant characteristic classes we have been considering.

\begin{prop} \label{p2.2}
We have
$$  \kappa_0 \Bigl ({\tdk(\symm)}{\tdki (\lieka \oplus \lieksa) }
\Bigr )
 = \td (\mred),   $$
where $\k_0 $ is the natural surjective ring homomorphism $\hk(\symm)
\to H^*(\mred) $.
Here, $\lieka $ denotes the product bundle $\symm \times \liek$
where $\liek$ is equipped with the adjoint action of $K$, and
$\lieksa $ denotes the  product bundle $\symm \times \lieks$
where $\lieks$ is equipped with the coadjoint action of $K$.
\end{prop}
\Proof The normal bundle $\nu(\zloc) $ to
$\mu^{-1}(0)$   (which is a submanifold of $\symm$ since
$0$ is a regular value
of $\mu$)  is isomorphic as an
equivariant bundle to $\lieksa$ (since $\mu: \symm \to
\lieks$ is an equivariant map).
Moreover, when $K$ acts freely on $\mu^{-1} (0)$, we have the following
decomposition
of $T\symm$ in terms of $K$-equivariant bundles:
$$ T \symm|_{\mu^{-1}(0)} = T (\mu^{-1} (0) ) \oplus \lieksa $$
and $  T (\mu^{-1} (0) ) = \proj^* T \mred \oplus \lieka $
where $\proj: \zloc  \to \mred$ is the natural projection.$\square$

The following is an immediate consequence of (\ref{4.9}):
\begin{lemma}
For $X \in \liet$ , the $T$-equivariant Todd class of $\lieka \oplus
\lieksa$ is given by
\beq \label{4.10}
\tdt(\lieka \oplus \lieksa)(X)  = \prod_{\g > 0 }
\frac{ \hg(X) ^2}{(1 - e^{ \isq \hg(X) } )(1 - e^{ - \isq \hg(X) } )},
\eeq
where the product is over the positive roots, and we have
introduced\footnote{The extra factors of $1/(2 \pi)$
are introduced because of the convention explained in Footnote 2
that weights $\b \in \liets$ satisfy
$ \b \in \Hom(\intlat, \ZZ) $ rather than
$\b \in \Hom(\intlat, 2 \pi \ZZ)$. Our definition of roots
is as in Lemma 3.1 of \cite{JK1}: thus the roots $\g$ satisfy
$\g(\intlat) \subset 2 \pi \ZZ$. In the terminology
of \cite{BD} (p. 185), the quantities $i \g: \liet \to i \RR$
are the {\em infinitesimal} roots whereas the corresponding
weights $\hg$ are the {\em real} roots.}
$\hg = \g/(2\pi)$.
\end{lemma}
When $K$ is abelian the bundle $\tdk(\lieka \oplus
\lieksa) $ is equivariantly
trivial as well as trivial, and so we have
in this case
\beq
  \kappa_0 (\tdk(\symm) )
 = \td (\mred).  \eeq

\renorm
\section{Review of the residue formula}
We now recall the main result (the residue formula, Theorem
8.1) of \cite{JK1}:
\begin{theorem} \label{t4.1}{[\cite{JK1}]} Let $\eta \in
\hk(\symm) $ induce $\eta_0 \in H^*(\mred)$. Then we have
\beq \label{4.1} \eta_0 e^{{}\omega_0} [\mred]
 = {\stabo C^K}  \res \Biggl (
\nusym^2 (X)
 \sum_{F \in \calf} \rfe(X) [d X] \Biggr ), \eeq
where $\stabo$ is the
order of the subgroup of $K$ that acts trivially on
$\zloc$, and
the constant $C^K$ is defined by
\beq \label{4.001} C^K = \frac{\isq^l}{(2 \pi)^{s-l}  |W| \vol(T)}. \eeq
 We have introduced $s = \dim K$ and
$l  = \dim T$. Also, $\calf$ denotes the set of components of the fixed point
set of $T$, and if $F$ is one of
these components then the meromorphic  function
$\rfe$ on $\liet \otimes \CC$ is defined by
\beq \label{4.01}  \onebl \rfe(X ) =  e^{i \mu_T(F) (X ) }
\int_F \frac{i_F^* (\eta(X ) e^{{}  \omega} )  }{e_F(X) }
 . \eeq
Here,  $i_F: F \to \symm$ is the inclusion
and $e_F$ is the $T$-equivariant Euler class of the normal
bundle to $F$ in $\symm$, which was defined at (\ref{3.021}).
 The polynomial
$\nusym: \liet \to \RR$ is defined by
$\nusym(X) = \prod_{\g > 0 } \g(X) $, where $\g$ runs
over the positive roots of $K$.
\end{theorem}
The general definition of the residue
$\res$ was given in Section 8 of \cite{JK1}. Here we shall
treat the case where $K$ has rank $1$, for which the results are as
follows.  See  Footnotes 2 and 3   for our conventions
on roots and weights.

\begin{corollary} \label{c4.2} {(\bf \cite{JK1}; \cite{Kalk},
\cite{wu})} In the
situation of Theorem \ref{t4.1}, let $K = U(1)$.
 Then
$$ \eta_0 e^{{}\omega_0} [\mred] =
{\isq\stabo}  \res_0 \Bigl (  \sum_{F \in \calf_+} \rfe(X)
d \lx \Bigr ) . $$
Here, the meromorphic  function $\rfe $ on  $ \CC$
 was defined
by (\ref{4.01}), and $\res_0 $ denotes the coefficient
of the meromorphic 1-form $d\lx/\lx$ on $\liek  \otimes \CC$,  where
$X \in \liek$ and
$\l$ is  the generator of  the weight
lattice of $U(1)$.    The set $\calf_+ $ is defined by
$\calf_+ = \{ F \in \calf: \mu_T(F) > 0 \}. $ The integer
$\stabo$ is as in Theorem \ref{t4.1}.
\end{corollary}

\begin{corollary} \label{c4.3}{\bf(cf. \cite{JK1}, Corollary 8.2)}
In the
situation of Theorem \ref{t4.1},
let $K = SU(2)$ or $K = SO(3) $. Then
$$ \eta_0 e^{{}\omega_0} [\mred] =
 \frac{\isq \stabo }{2} \res_0 \Bigl (  \hg(X)^2  \sum_{F \in \calf_+} \rfe(X)
d\l(X) \Bigr )  . $$
Here, $\res_0$, $\rfe$ and $\calf_+$ are as in Corollary
\ref{c4.2},  and  $\l   = \l_K \in
\liets$  is the generator of the weight lattice of $K$.
 We have $\l_{SO(3)}  = \hg$ and
$\l_{SU(2)} = \hg/2$,  where $\hg = \g/( 2 \pi)$ was defined
in terms of the positive root $\g$. The   integer
$\stabo$ is as in Theorem \ref{t4.1}.
\end{corollary}

We now specialize to the case $\eta_0 = \td (\mred)$.
Assume that $T$ acts at the
fixed point $F$ with weights $\bfj \in \liets$. From now on we
assume that the action of $K$ on $\zloc$ is
effective, so that $\stabo = 1$ in Theorem \ref{t4.1}.
\begin{prop} \label{p4.4} We have
\beq
\chtd  = {C^K} \res \Biggl ( \sum_{F \in \calf}
{\nusym^2(X)e^{\isq \mu_T(F)(X)}
 }{\tdti(\lieka \oplus \lieksa) (X) }
\eeq $$\times \int_F \frac{e^{{}\omega}
\tdt(\nu_F)(X) \td(F) }{e_F(X) }  \Biggr ). $$
This is equal to $\rr (\lred)$ provided $K$ acts
freely on $\zloc$.
\end{prop}
Here, the constant $C^K$ was  defined at (\ref{4.001}).
We have used the definitions of equivariant characteristic
classes given in Section 2. We have also decomposed
the restriction to $F$ of the $T$-equivariant Todd class of $\symm$ as
\beq \td_T (M) (X)   = \td_T (\nu_F)(X) \td (TF). \eeq
Here,  we
have used the multiplicativity of the Todd class and the fact that
the action of $T$ on $TF$ is trivial.
Then the Proposition follows immediately from Theorem
\ref{t4.1}.

The special case of Proposition \ref{p4.4} when $K = U(1)$ is:
\begin{prop} \label{p4.5} If $K =  U(1)$, we have
\beq \label{4.010} \chtd  =  {\isq} \res_0
\Bigl (   \sum_{F \in \calf_+} \expmf
\int_F \frac{e^{{}\omega}
\tdt(\nu_F)(X) \td(F) }{e_F(X) } d \lx   \Bigr ) \eeq
This is equal to $\rr (\lred)$ provided $K$ acts
freely on $\zloc$.
Here, $X \in \liet$ and $\res_0$ denotes the coefficient of
the meromorphic 1-form $d \lx/\lx$ on $\liet \otimes \CC$, where the element
$\l$ $ \in \liets$ is the generator of the weight lattice of $\liet$.
\end{prop}
The corresponding result for $K = SU(2)$ or $K = SO(3) $ is
\begin{prop} \label{p4.6} Let $K  = SU(2)$ or $K = SO(3)$. Then we
have
\beq \label{4.11} \chtd  =  \frac{\isq}{2} \res_0 \Bigl (
(1 - e^{\isq \hg(X) } ) (1 - e^{- \isq \hg(X) } )
\sum_{F \in \calf_+} \expmf\eeq
$$ \times \int_F \frac{e^{{}\omega}
\tdt(\nu_F )(X) \td(F)  }{e_F(X) }  d \l(X) \Bigr ).  $$
This is equal to $\rr (\lred)$ provided $K$ acts
freely on $\zloc$.
Here,  $ X \in \liet$, and
$\res_0$ denotes the coefficient of
the meromorphic 1-form $d \lx/\lx$ on $\liet \otimes \CC$, where the element
$\l  =
\l_K\in \liets$ is the generator of the weight
lattice of $K$;
also, $ \l_{SO(3) }  = $  $\hg = \g/(2 \pi)$ where $\g$ is the positive
root of $SO(3)$, and
$\l_{SU(2)} =
\hg/2$,  as in the statement of   Corollary \ref{c4.3}.
\end{prop}

Notice that it is valid to apply
the residue formula for groups of rank
one to formal equivariant cohomology classes
in this way, because both sides of the
formula send to zero all elements of $H^j_K(M)$
when $j>dim_{\bf R}(\mred)$.

\renorm
\section{The holomorphic Lefschetz formula}

We now describe the application of the holomorphic
Lefschetz theorem  in our situation. The  theorem  is proved by
Atiyah and Singer
(\cite{AS3}, Theorem 4.6), and   is based on results of Atiyah and
Segal \cite{AS2}: an exposition of the general result from which the
theorem follows is
given in
    Theorem 6.16 of \cite{BGV}.  A more general  equivariant index
theorem involving  equivariant cohomology is proved by
  Berline and Vergne in
\cite{BV2}. The following
statement  is in a form that will be convenient for us.
We introduce the notation that if $X \in \Lie(T)$
then
$t  = \exp (X) \in T$. For
 any
weight $\b$, we define   $t^\b$ as  $\exp ( 2 \pi i \b(X) ) $
$ \in U(1) \subset \CC^\times$, where the weights
$\b$ have been chosen to send the integer
lattice $\intlat $ in
$\liet$ to $ \ZZ \subset  \RR$.

\begin{theorem} \label{t5.0} {\bf (Holomorphic Lefschetz formula)}
Let $ t \in T $ be such that the
fixed point set of $t$ in $\symm$ is the same as the fixed point set
$\cup_{F \in \calf} F$ of $T$ in $M$;
then the character $\charr(t)$ of the virtual representation of $t$ on $\quant$
is given by
$$
\charr(t)  = \sum_{F \in \calf} \charr_F (t), $$
where
\beq \label{5.1} \charr_F (t) =
\int_F\frac{ i_F^* \cht(\call) (t) \td(F) }
{\prod_j (1 - t^{- \bfj} e^{ - c_1(\nfj) }) }
\eeq
$$ \onebl =
t^{\mf} \int_F \frac{e^{{}\omega} \td(F)  }
{\prod_j (1 - t^{- \bfj} e^{ - c_1(\nfj) }) }
 $$
Here,
the  $\bfj \in \Hom(T, U(1)) \subset \liets$ are the
weights of the action of $T$ on the normal bundle
$\nu_F$ of $F$ in $\symm$, and the $T$-moment map
$\mu_T$ is the composition of $\mu$ with restriction from
$\lieks$ to $\liets$.
\end{theorem}
\Proof Equation (\ref{5.1}) follows immediately from the statement
given in Theorem 4.6 of \cite{AS3}. We need only observe
 that the action of $t$ on the fibre of $\call$ above
any point in $F$ is given by multiplication by $t^{\mf}$.
Thus $i_F^* \ch_T (\call) (t) = e^{{}\omega} t^{\mf}. $$\square$

When the $T$ action has isolated fixed points, (\ref{5.1})
reduces to
\beq \label{5.2}
\charr_F (t) = \frac{t^{\mf} }
{\prod_j (1 - t^{- \bfj} ) } \eeq

In the general case, the structure of the right hand side of
(\ref{5.1}) is given as follows:
\begin{lemma} \label{l5.0} The expression
$$ \frac{1} {\prod_j (1 - t^{- \bfj} e^{ - c_1(\nfj) }) } $$
appearing in (\ref{5.1})
is given by
\beq \label{5.003}
\prod_{j} \sum_{r_j \ge 0 }
\frac{t^{- r_j \bfj} (e^{ - c_1(\nfj) } - 1 )^{r_j} }
{ (1 - t^{ - \bfj} )^{r_j  + 1} }. \eeq
In particular the only poles occur when $t^{\bfj}  = 1. $
\end{lemma}

\Proof This follows by examining for each $j$
$$ \frac{1}{ 1  - t^{ - \bfj}  e^{ - c_1 (\nfj) }   }
= \frac{1}{1 - y (1 +u)  }  =
\frac{1}{1 - y } \sum_{ r \ge 0 } \frac{y^r u^r }{ (1 - y)^r} $$
where $y = t^{ - \bfj} $
and $u = e^{ - c_1 (\nfj)}  - 1 $ is nilpotent. $\square$

We restrict from now on to the case $T = U(1)$, which is regarded as
embedded in $\hatc$ in the standard way.  We identify the
weights with integers by writing
them as multiples of the generator $\l$ of the weight
lattice of $U(1)$.

\begin{prop} \label{p5.1} The character $\charr(t )$ extends to
a holomorphic function on $\CC^\times = \hatc - \{ 0, \infty \}$.
 \end{prop}
\Proof This follows since $\charr$ is the character of a
finite dimensional  (virtual)
representation of $U(1)$, so it is of the
form $\charr(t) = \sum_{m\in \ZZ} c_m t^m $  for some
integer  coefficients $c_m$, finitely many of which are nonzero.
$\square$

The following is immediate:

\begin{prop} \label{p5.2} The expression $\charr_F$  given
in (\ref{5.1}) defines a meromorphic function
on $\hatc$  such that
$\sum_{F \in \calf} \charr_F(t) $  agrees with $\charr(t)$
on the open subset of $ U(1)$ consisting of those
$t$  whose action does not fix any point of
$\symm - \symm^T$.   Hence, by analyticity, $\charr(t) =
\sum_{F \in \calf} \charr_F(t)$ on an open set in $\hatc$ containing $
\CC^\times  - U(1)$.
\end{prop}

\begin{prop} \label{p5.3} The virtual dimension of the $T$-invariant
subspace of $\quant$ is given by
\beq \label{5.3}
\dim \quant^T =  \frac{1}{2 \pi \isq} \int_{|t| \in \G}
\frac{dt}{t} \sum_{F \in \calf} \charr_F(t), \eeq
where $\charr_F$ was defined after (\ref{5.1}).
Here, for any $\e > 0 $,
$\G = \{ t \in \hatc: |t| = 1 + \e  \}  \subset \Omega $
is a cycle in $\hatc $ on which   the $\charr_F$
have no poles.
\end{prop}
\Proof
This follows since
$$ \dim \quant^T = \frac{1}{2 \pi \isq} \int_{|t| = 1} \frac{dt}{t}
\charr(t) $$
$$ \bla =  \frac{1}{2 \pi \isq} \int_{|t|  \in \G} \frac{dt}{t}
\charr(t),  $$
and by applying Proposition \ref{p5.2} to identify $\charr $ with
$\sum_{F \in \calf } \charr_F $ on $\G$. $\square$

\noindent{\em Remark:} One obtains an equivalent formula
by defining $\G = \{ t \in \CC: |t| = 1 - \e \} $  for $0 < \e < 1$.

Let us now regard
\beq \label{5.03} h_F = \charr_F (t) \frac{dt}{t}  =
 \frac{dt}{t}
  t^{\mf} \int_F
\frac{ e^{{}\omega}
\td(F) }{\prod_j (1 - t^{ - \bfj}  e^{- c_1( \nfj) } ) }
\eeq
as a meromorphic 1-form on $\hatc$, whose  poles may occur
only at  $0$,
 $\infty $ and $\sss \in \wf$, where we define
\beq \label{5.06} \wf = \{ \sss \in U(1):
\sss^\bfj = 1  \mbox{ for some } \bfj \}. \eeq
 (This is true by inspection of
(\ref{5.2}) when the fixed point set of the action of $T$
consists of isolated points.
In the general case it follows from Lemma \ref{l5.0}.) The integral
(\ref{5.3}) then yields
\beq \label{5.4} \dim \quant^T  = -
\sum_{F \in \calf} \res_\infty h_F. \eeq

Let us examine the poles of $h_F $ on $\hatc$.
We have
\begin{lemma} \label{l5.4}
For a given $F$, let $n_{F,\pm}$
 be $\sum_{j: \pm \bfj > 0 } |\bfj|$.
If $\mf > - \nplf$ then
$\res_0  h_F = 0 $, while if $\mf < \nminf$ then
$\res_\infty h_F = 0 $. \end{lemma}
\Proof To study the residue at $0$, we assume
$|t| < 1 $, so that
$ (1 - t)^{-1} = \sum_{ n \ge 0 }  t^n $
and $(1 - t^{-1} )^{-1} = - t \sum_{n \ge 0 } t^n $.
For $r \ge 1 $ we examine
$$ \frac{ t^{ \mf} }{ \prod_j ( 1 - t^{ - \bfj} )^r }
\frac{dt}{t} $$
\beq  \label{5.5}
= t^{ \mf} (-1)^{l_+} t^{r \nplf}  \Bigl ( \prod_{j }
 \sum_{n_j \ge 0  } t^{|\bfj| n_j} \Bigr )^r \frac{dt}{t}, \eeq
where $l_+$ is the number of $\bfj$ that are positive.
It follows that if $\nplf + \mf > 0 $ then the residue at $0$
is zero.
A similar calculation yields the result for the residue at $\infty$.
$\square$

Recall that the action of $T$ on $M$ is said to be {\em quasi-free}
if it is free on the complement of the fixed point set of
$T$ in $M$.
The following is shown in \cite{DGP}:
\begin{lemma} \label{l5.3} The action of $T= U(1)$ on $\symm$
is quasi-free if and only if  the weights
are $\bfj = \pm 1$.
\end{lemma}

\begin{prop} \label{p5.5} If the action of $T$ is quasi-free, then
we have
\beq \label{5.7}
\res_\infty \sum_{F \in \calf} h_F =
 -  \sum_{F \in \calf_+ } \res_1 h_F . \eeq
Here,
 $\calf_+ = \{ F \in \calf: \mf > 0 \}$.
More generally the result is true if $\res_1 h_F $ is
replaced by
$\sum_{\sss \in \wf} \res_\sss h_F$, where the set
$\wf$ was defined at (\ref{5.06}).
\end{prop}

\Proof  Assume for simplicity that the action of
$T$ is quasi-free: the  proof of the general case is almost
identical. Lemma \ref{l5.4} establishes that
$$\res_\infty \sum_{F \in \calf} h_F = \sum_{F \in
\calf_+}
 \res_\infty h_F. $$
Also, if $F \in \calf_+$ then $\mf > - n_+$ so
$\res_0 h_F = 0 $; hence (\ref{5.7}) follows
because
the meromorphic 1-form $h_F$ has poles only
at $0, 1 $ and $\infty$ and their residues must sum to zero,
so $\res_1 h_F = - \res_\infty h_F $ when
$F \in \calf_+$.$\square$

\noindent{\em Remark:} Recall that $\mu_T(F)$ is never zero.

The following is an immediate consequence of combining Proposition
\ref{p5.5} with Proposition \ref{p5.3}:
\begin{corollary} \label{c5.7} If the action of $T = U(1)$ on
$\mred$ is quasi-free, we have $\rrt(\call) =
\sum_{F \in \calf_+} \res_1 h_F$.
More generally we have
$\rrt (\call)  = \sum_{F \in \calf_+} \sum_{\sss \in \wf}
\res_\sss h_F,  $
where $\wf$ was defined by (\ref{5.06}).
\end{corollary}

We now treat the cases $K  = SU(2)$ and
$K = SO(3)$.
 We shall first need the following
\begin{lemma} \label{l5.10}
There is no component $F$ of the fixed point set of $T$ on $\symm$
for which $\mf = 0 $.
\end{lemma}
\Proof Because the $K$ moment map is equivariant,
 $\mu(F)$ is fixed by
the action of $T$ on $\lieks$ for every $F \in \calf$;
 thus $\mu(F) \subset \liet$
(identifying $\liek$ with $\lieks$ and $\liet$ with
$\liets$ by the choice of an inner product), so that
$\mu(F)  = \mf$.  Thus $\mf  = 0 $ implies $\mu(F) = 0 $.
However, because $K$ acts locally freely on $\mu^{-1}(0)$,
no $F$ may intersect $\mu^{-1}(0)$.$\square$

We shall prove the
following result:
\begin{prop}  \label{p5.9}

\noindent{\bf (a)}  Suppose $M$ is connected,
and suppose $K = SO(3)$ acts on $\symm$ in such
a way that the action of $T$ is quasi-free.
Suppose also that there exists $F$ for which
$|\mu_T(F)| > 1$.
Then
$$ \rrk(\call)  = \half \sum_{F \in \calf_+}
\res_1 (2 - t - t^{-1}) h_F$$ where the meromorphic 1-form
$h_F $ on $\hatc$ was defined by (\ref{5.03}).
More generally we have when $K = SO(3)$ (provided there
exists $F$ for which
$|\mu_T(F)| > 1$) that
$$ \rrk(\call) =
\half \sum_{F \in \calf_+} \sum_{\sss \in \wf}
\res_\sss  (2 - t - t^{-1}) h_F,$$
where $\wf$ was defined by (\ref{5.06}).

\noindent{\bf (b)} Let
$K = SU(2)$, and suppose that there is an $F$ for which
$|\mu_T(F) | > 2, $ and also that there is no $F$ with
either $\mu_T(F) = 1$ and $\nplf = 1$ or
$\mu_T(F) = - 1 $ and $\nminf = 1$.
Then we have that

$$ \rrk(\call) =
\half \sum_{F \in \calf_+} \sum_{\sss \in \wf}
\res_\sss  (2 - t^2 - t^{-2}) h_F.$$
\end{prop}

\Proof
{\bf (a)}:If $K = SO(3)$, we have by the Weyl integral formula for Lie groups
that
\beq \label{5.07}\dim \quant^K = \frac{1}{\vol K}
\int_{k \in K} dk \charr(k) = \frac{1}{|W|}  \frac{1}{2 \pi \isq}
\int_{t \in T} \frac{dt}{t} (1 - t) (1 - t^{-1}) \charr(t) \eeq
$$ \bla = \half \frac{1}{2 \pi \isq} \int_{t \in \G} \frac{dt}{t}
(2 - t - t^{-1}) \charr(t) $$
$$ \bla = \half \frac{1}{2 \pi \isq} \int_{t \in \G} \frac{dt}{t}
(2 - t - t^{-1}) \sum_{F \in \calf} t^{\mf} \int_F
\frac{ e^{{}\omega}
\td(F) }{\prod_j (1 - t^{ - \bfj}  e^{- c_1( \nfj)} ) }. $$
(In (\ref{5.07}), the factor $(1 - t) (1 - t^{-1})$ is the
volume of the conjugacy class  of $K$ containing $t$, in a
normalization
where $\vol K = 1$:  see for
instance \cite{BD}, (IV.1.11).)
By the previous argument (Lemmas
 \ref{l5.0} to \ref{l5.4}), this
is equal to
\beq \label{5.8}
\dim \quant^K = - \half
\sum_{F \in \calf} \res_\infty (2- t - t^{-1}) h_F \eeq
 where the meromorphic 1-form
$h_F $ was defined at (\ref{5.03}).
By the proof of Lemma \ref{l5.4} this becomes
\beq \label{5.9}
\half \bigl ( 2 \res_1 \sum_{F \in \calf_+} h_F
- \res_1 \sum_{F \in \calf, \mf + 1 \ge \nminf} t h_F
- \res_1 \sum_{F \in \calf, \mf - 1 \ge \nminf} t^{-1} h_F \bigr ). \eeq
By Lemma \ref{l5.10},
$\mf \ne 0 $ for any $F$; thus it suffices  to check that
\beq \label{5.91}
\sum_{F \in \calf: \onebl \mu_T(F) + 1 \ge \nminf} \res_1 (t h_F)
 =
\sum_{F \in \calf: \onebl \mu_T(F) > 0 } \res_1 (t h_F)  \eeq
and likewise that
\beq \label{5.92}
\sum_{F \in \calf: \onebl \mu_T(F) -  1 \ge \nminf} \res_1 (t^{-1} h_F)
 =
\sum_{F \in \calf: \onebl \mu_T(F) > 0 } \res_1 (t^{-1} h_F).  \eeq
Equation (\ref{5.91})  follows by the proof of
 Lemma \ref{l5.4}, unless $\nminf = 0 $ and
$\mf = - 1$: we  apply the proof, replacing $\mf$ by
$\mf + 1 $ and using the fact that  for any $r \in \ZZ$,
\beq \label{5.94} \res_0(t^r h_F) =  \res_\infty
(t^r h_F) = 0  ~~~\mbox{if}~~~\mf+ r \in
[ - \nplf + 1, \nminf - 1]. \eeq
Likewise, equation (\ref{5.92})  follows unless $\nplf = 0$ and $\mf = 1$.
However
if $(\mf,\nminf)  = (-1,0)$
then $F$ gives a local minimum of
$\mu_T$.
Since $M$ is connected and $\mu_T$ is a perfect
Morse function (\cite{Ki1}, (5.8)), the local minimum must
be a global minimum,
contradicting the assumption that there exists an $F'$ for which
$|\mu_T(F')| > 1$.\footnote{Recall that by Weyl symmetry
there exists $F'$ for which $\mu_T (F') > 1 $ if and only if there exists
$F'$
 for which $\mu_T (F') < - 1 $.}
 Similarly the case $(\mf, \nminf) = (1,0) $  gives
a     maximum of $\mu_T$ and
hence cannot occur.

\noindent{\bf (b)} If $K = SU(2)$ we obtain instead of
(\ref{5.9})
\beq \label{5.93} \dim \quant^K =
\half \bigl ( 2 \res_1 \sum_{F \in \calf_+} h_F
- \res_1 \sum_{F \in \calf, \mf + 2\ge \nminf} t^2 h_F
- \res_1 \sum_{F \in \calf, \mf - 2  \ge \nminf} t^{-2} h_F \bigr ). \eeq
Using (\ref{5.94}) we find that  the second sum in
(\ref{5.93}) is equal to
$- \res_1 \sum_{F \in \calf, \mf > 0 } t^2 h_F $ except
when $\mf,\nminf)$ is $(-2,0)$, $(-1,0)$ or $(-1,1)$.
Likewise the third sum is equal to
$- \res_1 \sum_{F \in \calf, \mf > 0 } t^{-2} h_F $ except when
$(\mf, \nplf)$ is $ (2,0)$, $(1,0)$ or $(1,1)$.
The first, second, fourth and fifth of these six cases are excluded
if we assume that there is some $F$ with $|\mf| > 2$. $\square$

\noindent{\em Remark:} The technical hypothesis in
Proposition \ref{p5.9}(a)  that there should exist
$F$ for which $|\mf| \ne 1$ can be satisfied
by replacing $\call$ by $\call^k$
with $k \ge 2$. Similarly the technical hypotheses in
\ref{p5.9}(b) can be satisfied by taking $k \ge 3.$

\noindent
\renorm
\section{Identification with the residue formula}

In this section we shall  assume the weights are $\bfj = \pm 1$, so that
the action of $T$ on $M$ is quasi-free.
In order to treat the general case one needs to use
Kawasaki's
Riemann-Roch theorem for orbifolds \cite{Kaw}: we do this
in the next section.

Let us examine the residue $ \res_1 h_F$ in the case $K =  U(1)$.
We denote a generator of the weight lattice of $\liet$ by
$\l$, and replace the
parameter $t$ (in a small neighbourhood of
$1 \in \hatc$)  by
\beq \label{6.001} t = e^{i \l(X) }  \eeq
(where  $X \in \liet \otimes
\CC$
is in a small neighbourhood of $0$ in $\liet \otimes \CC$),
so that
$$ \frac{ dt}{t}  =  i { d \lx} $$
defines  a meromorphic 1-form on $\liet \otimes \CC$.
(The substitution (\ref{6.001}) differs from the
substitution used in Section 4, where we set
$t = e^{2 \pi i \l(X)}$: however the value of the residue
obviously is independent of  which of these substitutions
is used, and the substitution (\ref{6.001})  yields the formulas in
Section 3.)

We then find that
\beq \res_1 h_F = \isq \res_0
\Bigl (    \expmf
\int_F \frac{e^{{}\omega} \td(F)} {\prod_{j } (1 - e^{- \isq
\bfj( X) - c_1 (\nfj) } ) } d  \lx \Bigr ) \eeq
\beq \label{6.1}
\bla = \isq \res_0  \Bigl (
\expmf \int_F
\frac{  e^\omega \tdt (\nu_F )(X)\td(F)   }
{ \prod_j (\isq \bfj(X) + c_1 (\nfj)  ) }
d \lx \Bigr ),\eeq
\beq \label{6.2} \bla =
\isq  \res_0  \Bigl (
 e^{ \isq \mf(X) } \int_F  \frac{ e^{{} \omega}
\tdt (\nu_F)(X)  \td(F) }{e_F(X) } d \lx
\Bigr )  \eeq
where $\res_0 $ denotes the coefficient of $\srlx$.
Combining (\ref{6.2}) with
Proposition  \ref{p4.5},
  one obtains
\begin{prop} \label{p6.1}  We have
\beq \chtd = \sum_{F \in \calf_+}
\res_1 h_F. \eeq
This equals $\rr(\lred) $ provided $K = U(1) $ acts freely on $\zloc$.
 \end{prop}
Comparing  Proposition \ref{p6.1} with Corollary \ref{c5.7}, we have
\begin{prop}  \label{p6.1p}
Let the action of $K  =  U(1) $  on $\symm$ be quasi-free  (which
implies $K$ acts freely on $\mu^{-1}(0) $). Then
$\rrk(\call) = \rr(\lred)$. \end{prop}

To treat $K = SO(3)$ and $K=SU(2)$,
using the  substitution (\ref{6.001}) in $\res_1 (2 - t - t^{-1}) h_F$,
we recover the
right hand side of Proposition \ref{p4.6}:

\begin{prop} \label{p6.04} Let $K = SO(3)$ or $SU(2)$ act on
$\symm$.
Then
$$\chtd  = \half \sumpl \res_1 (2 - t - t^{-1}) h_F. $$
This equals $\rr(\lred) $ provided $K$ acts freely on $\zloc$.
\end{prop}
Combining Proposition \ref{p6.04} with  Proposition \ref{p5.9}
we get\footnote{We do not treat
$K = SU(2)$, since in this case the
 action of $T$ can only be quasi-free  if all the $F$ are in
the fixed point set of $K$.
For unless
$F$ is fixed by all of $K$,
 the orthocomplement $\liet^\perp$ of $\liet$ (equipped with the
adjoint action) injects into the
normal bundle $\nf$ under the action of $K$.
There is thus a subbundle  of $\nf$  on which $T$ acts with
weight $2$ or $-2$, and so the action of $T$ cannot be quasi-free
by Lemma \ref{l5.3}.}
\begin{prop} \label{p6.04p} Let $K = SO(3) $ act on $\symm$
 in such a way that the action of $T$ is quasi-free
(which
implies $K$ acts freely on $\mu^{-1}(0) $). Then
$\rrk(\call) = \rr(\lred)$. \end{prop}

Thus we have
\begin{theorem} \label{t6.2} Suppose $T = U(1) $ and
either $K = U(1) $ or $K = SO(3)$. Suppose
$K$ acts in a Hamiltonian fashion on the  K\"ahler
manifold\footnote{{\rm As described in the Introduction, it actually suffices
to assume that $M$ is equipped with a  $K$-invariant
 {\em almost} complex structure compatible with
the symplectic structure. This applies likewise to Theorem
\ref{t7.2} below.}}
$\symm$,  in such a way that the action of $T$ is
quasi-free. If $K = SO(3)$, suppose also that there exists $F$
for which $|\mf| > 1$.
We assume a moment map $\mu$ for the action of $K$
has been chosen in such a way that $0$ is a regular value of $\mu$.
Then $\rrk(\call)  = \rr(\lred).$
\end{theorem}

\section{Kawasaki's Riemann-Roch theorem}
In this final section we sketch the proof of the Guillemin-Sternberg
result
$\rr(\lred) = \rrk(\call)$ when $K$ has rank one,  without the assumption
that the action of $T$ is quasi-free. In this more general case,
$\mred$
is an orbifold and $\lred$ an orbifold bundle. The Riemann-Roch
number of $\lred$ is then given by applying Kawasaki's Riemann-Roch
theorem for orbifolds.
We state Kawasaki's result only as it applies
in our particular situation:
the special case when $K=T$ in fact appears in earlier
work of Atiyah.
\begin{theorem} \label{t7.1}  {(\bf  Atiyah \cite{A:ell}; Kawasaki
\cite{Kaw})} The Riemann-Roch
number of the orbifold bundle $\lred$ is given by
\beq \label{7.1} \rr(\lred) =
\chtd + \sum_{ 1\neq\sss \in \cals}
 \sum_{a \in \cala_s} \frac{1}{n_{s,a} }
\int_{\msoab} \isa. \eeq
Here, $\cals$ is a set
of representatives $s \in T $ for the conjugacy classes in
$K$ of elements
whose fixed point set $\ms$ is strictly larger than the fixed point
set of any subgroup of $K$ of dimension at least one.
 The components
of $\ms$ are denoted $\msa$, where $a \in \cala_s$;
we introduce $\msoa = \msa \cap \zloc$,
and $\msoab = \msoa/K_s$ where $K_s$ is the centralizer of $s$ in $K$.
The positive integer $n_{s,a} $ is the order of the
stabilizer of the action of $K_s$ at a generic point
of $\msa$.
 The class $\isa $ $ \in H^*(\msoab) $
is defined by
\beq \label{7.01} \isa = \frac{\ch({\call}^a_{s,red}) s^{\mu_a} \td(\msoab) }
{\prod_{k \in \kappind}  \bigl ( 1 - s^{ - \betak} e^{- \cnus}  \bigr
) }.   \eeq
Here, $\mu_a$ is the weight of the action of $s$ on the
fibre of $\call$ over any point in
$\msa$ and ${\call}^a_{s,red}$ is the induced orbifold bundle on $\msoab$.
If $\nu(\msoab)$ denotes the orbifold bundle
which is the pullback to $\msoab$ of the normal
to the image of the natural map from $\msoa$ to $\zloc$, we
decompose $\nu(\msoab)$ as a formal sum of line
bundles
\beq
\nu(\msa) = \oplus_{k \in \kappind}  \nusak , \eeq
and denote by $\betak \in \ZZ$ the weight of the action
of $s$  on the formal line
subbundle of the normal bundle to $\msoa$ in $\zloc$
corresponding to $\nusak$.
\end{theorem}
We can use this Theorem to prove Guillemin and Sternberg's
result for groups of
rank one, by identifying the additional  terms on the right hand side of
(\ref{7.1}) with the additional residues at the points
$1 \ne $ $\sss \in\wf$
that appear in the statement of Proposition
\ref{p5.5}  when the action of $T$ is not quasi-free.
Meinrenken uses Kawasaki's theorem in a different way to eliminate the
 quasi-free action  hypothesis from the proof given by
Guillemin in \cite{G}: see \cite{mein}, Remark 1 following Theorem
2.1.

The proofs of Corollary \ref{c5.7} and Proposition 4.11
give when $K=T$
\beq \label{7.2} \rrk(\call) =
\sum_{s \in \cals} \sum_{a \in \cala_s}
\bigl ( \sumfms   \res_s h_F \bigr ), \eeq
and when $K=SO(3)$ or $SU(2)$
\beq \label{7.25} \rrk(\call) = \sum_{s \in W\cals} \sum_{a \in \cala_s}
\bigl( \sumfms \res_s (2 - t - t^{-1}) h_F/2 \bigr ), \eeq
where the meromorphic 1-form $h_F$ on $\hatc$ was defined
at (\ref{5.03}) and $W$ is the Weyl group of $K$. The terms in the second
sum indexed by different elements $s$ of the same Weyl group orbit are equal,
so the sum can be rewritten as a sum over  $\cals$ instead of $W\cals$.
We know from Proposition \ref{p6.1} (a consequence
of applying the residue formula  to the class $\ch(\lred) \td(\mred)$
on $\mred$) that the term in each of these sums indexed by $s=1$ is
\beq \label{7.3} \chtd . \eeq
To deal with the other terms we apply the residue formula (Theorem
\ref{t4.1}) to the action of
$K_s$ on the symplectic manifold $\msa$: in the notation
of that theorem, we  introduce an
appropriate equivariant cohomology class $\eta e^{\bom
 }  = $
 ${\cali}^{s,a}_{K_s} \in H^*_{K_s}(\msa)$ which
descends on the symplectic quotient  $\msoab$ to
$\isa$ $ =
\eta_0 e^{\omega_0} $.
 When $K=T$ the class $\isat$ is  given by
\beq \label{7.03}  \isat = {\cht(L) s^{\mu_a} \tdt(\msa) }
\Bigl  ( \prod_{k \in \kappind}  \bigl ( 1 - s^{ - \betak} e^{- \cnust}  \bigr
)^{-1}     \Bigr ) \eeq
where $\cnust$ is the $T$-equivariant
first Chern class of the
virtual line bundle $\nusak$. In the
other cases ${\cali}^{s,a}_{K_s}$ is defined
similarly, using Proposition 2.2 applied to $K_s$.
 This yields for each $s \in \cals $ and
$a \in \cala_s$ that the
term in the right hand side of  (\ref{7.2}) or (\ref{7.25})
indexed by $s$ and $a$ is\beq \label{7.4}
\frac{1}{n_{s,a}} \int_{\msoab} \isa. \eeq
(Here, the factor $n_{s,a}$ is the order of the subgroup of
 $K_s$ that acts trivially on  $\msa$: see the statement of
Theorem \ref{t4.1}.)
Substituting (\ref{7.4}) in (\ref{7.2}) or (\ref{7.25}) we recover the
right hand side of (\ref{7.1}). Thus we  obtain  the
Guillemin-Sternberg result in the special case when $K $ has rank one:
\begin{theorem} \label{t7.2} Suppose
that a compact group $K$ with maximal torus
$T = U(1) $
 acts in a Hamiltonian fashion on the  K\"ahler
manifold $\symm$,
in such a way that $0$ is a regular value of $\mu$.
Then if the hypotheses of Proposition 4.11(a) and (b) are satisfied
$$\rrk(\call)  = \rr(\lred).$$
\end{theorem}

\noindent{\bf Remark:} It has been
pointed out to us by M. Vergne that there are examples
(such as the action of $SU(2)$ on
the complex projective line) to show that this result
is not true without some hypotheses such as those of Proposition 4.11.

\end{document}